\documentclass[11pt]{article}
\usepackage{amsfonts,amssymb, graphicx,a4,bm}
\usepackage[english]{babel}
\usepackage{amsmath,amsthm,epsfig}
\def\theequation{\arabic{equation}}

\newtheorem{theorem}{Theorem}
\newtheorem{proposition}[theorem]{Proposition}

\begin{document}

\def\reff#1{(\protect\ref{#1})}

\let\a=\alpha \let\b=\beta \let\ch=\chi \let\d=\delta \let\e=\varepsilon
\let\f=\varphi \let\g=\gamma \let\h=\eta    \let\k=\kappa \let\l=\lambda
\let\m=\mu \let\n=\nu \let\o=\omega    \let\p=\pi \let\ph=\varphi
\let\r=\rho \let\s=\sigma \let\t=\tau \let\th=\vartheta
\let\y=\upsilon \let\x=\xi \let\z=\zeta
\let\D=\Delta \let\F=\Phi \let\G=\Gamma \let\L=\Lambda \let\Th=\Theta
\let\O=\Omega \let\P=\Pi \let\Ps=\Psi \let\Si=\Sigma \let\X=\Xi
\let\Y=\Upsilon

\global\newcount\numsec\global\newcount\numfor
\gdef\profonditastruttura{\dp\strutbox}
\def\senondefinito#1{\expandafter\ifx\csname#1\endcsname\relax}
\def\SIA #1,#2,#3 {\senondefinito{#1#2}
\expandafter\xdef\csname #1#2\endcsname{#3} \else
\write16{???? il simbolo #2 e' gia' stato definito !!!!} \fi}
\def\etichetta(#1){(\veroparagrafo.\veraformula)
\SIA e,#1,(\veroparagrafo.\veraformula)
 \global\advance\numfor by 1
 \write16{ EQ \equ(#1) ha simbolo #1 }}
\def\etichettaa(#1){(A\veroparagrafo.\veraformula)
 \SIA e,#1,(A\veroparagrafo.\veraformula)
 \global\advance\numfor by 1\write16{ EQ \equ(#1) ha simbolo #1 }}
\def\BOZZA{\def\alato(##1){
 {\vtop to \profonditastruttura{\baselineskip
 \profonditastruttura\vss
 \rlap{\kern-\hsize\kern-1.2truecm{$\scriptstyle##1$}}}}}}
\def\alato(#1){}
\def\veroparagrafo{\number\numsec}\def\veraformula{\number\numfor}
\def\Eq(#1){\eqno{\etichetta(#1)\alato(#1)}}
\def\eq(#1){\etichetta(#1)\alato(#1)}
\def\Eqa(#1){\eqno{\etichettaa(#1)\alato(#1)}}
\def\eqa(#1){\etichettaa(#1)\alato(#1)}
\def\equ(#1){\senondefinito{e#1}$\clubsuit$#1\else\csname e#1\endcsname\fi}
\let\EQ=\Eq
\def\0{\emptyset}

\def\pp{{\bm p}}\def\pt{{\tilde{\bm p}}}


\def\\{\noindent}
\let\io=\infty

\def\VU{{\mathbb{V}}}
\def\EE{{\mathbb{E}}}
\def\GI{{\mathbb{G}}}
\def\TT{{\mathbb{T}}}
\def\C{\mathbb{C}}
\def\CC{{\mathcal C}}
\def\II{{\mathcal I}}
\def\LL{{\cal L}}
\def\RR{{\cal R}}
\def\SS{{\cal S}}
\def\NN{{\cal N}}
\def\HH{{\cal H}}
\def\GG{{\cal G}}
\def\PP{{\cal P}}
\def\AA{{\cal A}}
\def\BB{{\cal B}}
\def\FF{{\cal F}}
\def\v{\vskip.1cm}
\def\vv{\vskip.2cm}
\def\gt{{\tilde\g}}
\def\E{{\mathcal E} }
\def\I{{\rm I}}
\def\rfp{R^{*}}
\def\rd{R^{^{_{\rm D}}}}
\def\ffp{\varphi^{*}}
\def\ffpt{\widetilde\varphi^{*}}
\def\fd{\varphi^{^{_{\rm D}}}}
\def\fdt{\widetilde\varphi^{^{_{\rm D}}}}
\def\pfp{\Pi^{*}}
\def\pd{\Pi^{^{_{\rm D}}}}
\def\pbfp{\Pi^{*}}
\def\fbfp{{\bm\varphi}^{*}}
\def\fbd{{\bm\varphi}^{^{_{\rm D}}}}
\def\rfpt{{\widetilde R}^{*}}

\def\tende#1{\vtop{\ialign{##\crcr\rightarrowfill\crcr
              \noalign{\kern-1pt\nointerlineskip}
              \hskip3.pt${\scriptstyle #1}$\hskip3.pt\crcr}}}
\def\otto{{\kern-1.truept\leftarrow\kern-5.truept\to\kern-1.truept}}
\def\arm{{}}
\font\bigfnt=cmbx10 scaled\magstep1

\newcommand{\card}[1]{\left|#1\right|}
\newcommand{\und}[1]{\underline{#1}}
\def\1{\rlap{\mbox{\small\rm 1}}\kern.15em 1}
\def\ind#1{\1_{\{#1\}}}
\def\bydef{:=}
\def\defby{=:}
\def\buildd#1#2{\mathrel{\mathop{\kern 0pt#1}\limits_{#2}}}
\def\card#1{\left|#1\right|}
\def\proof{\noindent{\bf Proof. }}
\def\qed{ \square}
\def\trp{\mathbb{T}}
\def\trt{\mathcal{T}}
\def\Z{\mathbb{Z}}
\def\be{\begin{equation}}
\def\ee{\end{equation}}
\def\bea{\begin{eqnarray}}
\def\eea{\end{eqnarray}}

\title {An Improvement of the Lov\'{a}sz Local Lemma via
Cluster Expansion}

\author{Rodrigo Bissacot$^{1,2}$, Roberto Fern\'{a}ndez$^{2,3}$,  Aldo Procacci$^1$ and Benedetto
Scoppola$^4$ \\\\
\footnotesize{$^1$Dep. Matem\'atica-ICEx, UFMG, CP 702
Belo Horizonte - MG, 30161-970 Brazil}
\\
\footnotesize{$^2$Labo. de Maths Raphael Salem, Universit\'e de Rouen, 76801,  France}\\
\footnotesize{$^3$Department of Mathematics, Utrecht University, P.O. Box 80010 3508 TA Utrecht\footnote{On leave from
Laboratoire de Math\'ematiques Raphael SALEM,
UMR 6085 CNRS-Universit\'e de Rouen, France}}\\
\footnotesize{$^4$Dipartimento di Matematica - Universita Tor Vergata di Roma, 00133 Roma, Italy}\\
\tiny{emails: {rodrigo.bissacot@gmail.com};~ {R.Fernandez1@uu.nl};}\\
\tiny{
{ aldo@mat.ufmg.br};~
{scoppola@mat.uniroma2.it}}
}
\date{}
\maketitle

\numberwithin{equation}{section}

\begin{abstract}
An old result by Shearer relates the Lov\'asz Local Lemma with  the
independent set polynomial on graphs, and consequently, as observed by Scott and Sokal, with the partition function of the  hard core
lattice gas on graphs. We use this connection and a recent
result on the analyticity of the logarithm of the partition function
of the abstract polymer gas to get  an improved version of  the Lov\'asz Local Lemma.
As applications we obtain tighter bounds on conditions for the existence of latin transversal matrices and the satisfiability of $k$-SAT forms.
\end{abstract}

\section{Introduction}

The main aim of this paper is to outline how techniques developed within statistical mechanics
can be applied to improve combinatorial results proved by means of the Lovasz Local Lemma (LLL).
These improvements rely on three major contributions:
\begin{itemize}
\item[(i)] Shearer's relation \cite{Sh} between the LLL and independent-set
polynomials;

\item[(ii)] Scott and Sokal's subsequent connection \cite{So} with lattice-gas
partition functions, and

\item[(iii)] relatively recent results in \cite{FP} on the convergence radius of logs
of partition functions.
\end{itemize}

While point (iii) contains the main mathematical tool, the steps needed to exploit its combinatoric consequences
by the way of the LLL, are far from obvious.  Our note is intended,
then, to fulfill a double role: First, to offer the
stochastic-combinatoric community an improvement of the LLL whose
power is illustrated in well known examples.  Second, to present a
simple and explicit road-map showing how improvements in cluster
expansion estimations lead to improvements in the LLL.  This can be seen as a streamlined presentation of the
connections exploited in \cite{So}.
Thus, the results of our paper yield, on the one hand, an effective tool to solve some combinatoric problems,
and, on the other hand, an illustration of how
statistical mechanical methods can be profitably used in seemingly
unrelated ---not even probilistic--- problems.

This short note is organized as follows.  In the next section we recall the
connections between LLL and statistical mechanics, and we state the improvement of LLL
based on \cite{FP}. Applications are presented in the following section. Some technical aspects of the basic result of this paper
are presented in appendix.
\section{Results}

\hspace{0.7cm}

One of the more powerful tools used in the probabilistic method in combinatory
is the so-called {\it Lovasz local Lemma}, first proved by Erd\"{o}s and Lov\'{a}sz in \cite{ES}.
To state this lemma we need several preliminary definitions.

Hereafter $|U|$ denotes the cardinality of a finite set  $U$.
Let $X$ be a finite set and $\{A_x\}_{x\in X}$ be a family of events on some probability space,
each of which having probability $\mathbb P(A_x)=p_x$ to occur.
A graph $G$ with vertex set $V(G)=X$
is a {\it dependency graph} for the family of events
$\{A_x\}_{x\in X}$  if, for each $x\in X$,
$A_x$ is independent of all the events in the $\sigma$-algebra
generated by $\{ A_y: y \in
{X}\backslash \Gamma^{*}_G(x) \}$, where $\Gamma_G(x)$ denotes the
vertices of $G$ adjacent to $x$ and
$\Gamma^{*}_G(x)=\Gamma_G(x)\cup\{x\}$. For any  non empty $S\subset X$ we also denote by  $\Gamma_G(S)$ the neighborhood of $S$, i.e.,
 $\Gamma_G(S)=\{x\in X\backslash S:\mbox{$x$ adjacent to $y$ in $G$ for some $y\in S$}\}$.

Let now $\bar A$ be the event $\bar A=\cap_{x\in X}\bar A_x$,
where  $\bar A_x$ is the complement event of $A_x$.
In words, $\bar A$ is the event that none of the events
$\{A_i\}_{i\in {X}}$  occurs. The Lovaz local lemma gives a sufficient
criterion to guarantee that $\bar A$ has positive probability (and hence it exists).


\begin{theorem}[{\bf Lov\'asz local lemma}]
\label{LLLold} Suppose that $G$ is a dependence
graph for the family of events $\{A_x\}_{x\in X}$   with probability
$\mathbb{P}(A_x)=p_x$ and there exist $\{\m_x\}_{x\in X}$ real
numbers in $[0,+\infty)$ such that, for each $x\in X$,
\begin{equation}\label{lov}
p_x \;\leq\; R_x\;=\; {\mu_x \over \varphi_x(\bm \mu)}
\end{equation}
with\
\begin{eqnarray}
\varphi_x(\bm \mu) & =& \prod_{y\in \Gamma^*_G(x)} (1+\m_y)\\ \label{aab}
& = & \sum_{R\subset  \Gamma^*_H(x)} \prod_{x\in R}\mu_x  \label{ab}
\end{eqnarray}
 Then
 \begin{equation}\label{a}
\mathbb{P}(\cap_{x\in X}\bar A_x)\;\ge\; \prod_{x\in X} {1\over (1+\m_x)}\,> \,0
\end{equation}
\end{theorem}

Here and in the sequel the sum    over all subsets  in the  r.h.s. of \reff{ab} includes the empty set and we agree that $\prod_{x\in \emptyset}\mu_x\equiv 1$
\vskip.2cm
\\{\bf Remark}. In the literature the Lov\'asz local lemma above is usually written in terms of variables $r_x= {\m_x/( 1+\m_x)}\in [0,1)$,
so that condition (\ref{lov}) is replaced by  $p_x \;\leq\; \; {r_x \prod_{y\in \G_G(x)}( 1-r_y)}$ and bound \reff{a} is replaced by
$\mathbb{P}(\cap_{x\in X}\bar A_x\ge\prod_{x\in X}(1-r_x)$.
The present formulation, however, shows more directly the improvements contained in Theorem \ref{LLL}.

\vskip.2cm

The connection between the
Lovasz local lemma and the statistical mechanics of the hard-core gas
has been pointed out by Scott and Sokal in \cite{So}.  They used and clarified an
old result by Shearer \cite{Sh} which can be stated as follows.

\begin{theorem}[{\bf Shearer}]\label{Sh}
Suppose that $G$ is a dependence
graph for the family of events $\{A_x\}_{x\in X}$ each one having probability
$\mathbb{P}(A_x)=p_x$.
Let $S\subset X$.
Define
\begin{equation}
P(S)=\sum_{U:  ~S\subset U \subset X\atop U\ {\rm indep\ in}\ G}
\!\!\!(-1)^{|U|-|S|}\prod_{x\in U}p_x
\end{equation}
If $P(S)\ge0$ for all $S\subset X$ then
\begin{equation}
\label{eq:r1}
\mathbb P(\bar A)\;\ge\; P(\emptyset)\;\equiv
\sum_{U:   U \subset X\atop U\ {\rm indep\ in}\ G}
\!\!\!(-1)^{|U|}\prod_{x\in U}p_x'.
\end{equation}
\end{theorem}

\vv
Theorem \ref{Sh} can be rewritten in terms of the statistical mechanical {\it partition function} $Z_G(\bm w)$ of the hard-core lattice gas with complex activities
${\bm w}=\{w_x\}_{x\in {X}}$ (with $w_x\in \mathbb{C}$ for all $x\in {X}$) on the graph $G$,
defined by
\begin{equation}\label{part}
Z_G(\bm w)=\sum_{R\subset  {X}\atop R\ {\rm indep\ in}\ G}
\prod_{i\in R}w_x
\end{equation}

It is now  immediate  to see that
\be\label{p0}
P(\emptyset)=Z_G(-\bm p)
\ee
where ${ -\bm p}=\{-p_x\}_{x\in {X}}$.
Moreover, for any non-empty $S\subset X$,
\begin{equation}\label{ps}
P(S)=\sum_{U:  ~S\subset U \subset {X}\atop U\ {\rm indep\ in}\ G}
\!\!\!(-1)^{|U|-|S|}\prod_{x\in U}p_x=
\sum_{R\subseteq {X}\backslash (S\cup\Gamma_G(S)) \atop R~{\rm indep\ in}\ G}\prod_{x\in R}(-p_x)\prod_{y\in S}p_y
=Z_G(-\bm p_S)\prod_{y\in S}p_y
\end{equation}
where ${ -\bm p_S}=\{-p^S_x\}_{x\in {X}}$ and
\begin{equation}
p^S_x=
\begin{cases} 0 & \text{$x\in S\cup\Gamma_G(S)$}
\\
p_x &\text{otherwise}
\end{cases}
\end{equation}






Formulae (\ref{p0}) and (\ref{ps}) and Theorem 2 imply the following proposition.

\begin{proposition}\label{elo}. Suppose that $G$ is a dependence
graph for the family of events $\{A_x\}_{x\in X}$ each one having probability
$\mathbb{P}(A_x)=p_x$.
Let $Z_G(\bm w)$ be the partition function of of the hard-core gas on $G$ with complex activities
${\bm w}=\{w_x\}_{x\in {X}}$ (with $w_x\in \mathbb{C}$ for all $x\in {X}$).
Suppose that $\log Z_G(\bm w)$ is an absolutely convergent series
in the polydisc $|w_x|\le p_x$,  Then
\begin{equation}
\label{eq:r2}
\mathbb P(\bar A)\;\ge\; Z_G(-\bm p)\;>\; 0\;.
\end{equation}
\end{proposition}

\\{\bf Proof}.  The restriction of  $Z_G(\bm w)$  to  $K_{\bm p}=\prod_{x\in X}[-p_x,p_x]$
is a real polynomial function  and it is positive in  $K^+_{\bm p}=\prod_{x\in X}[0,p_x]$.
Moreover, since $\log Z_G(\bm w)$ is an absolutely convergent series in the polydisc $|w_x|\le R_x$, the partition function $Z_G(\bm w)$ has no zeroes in the same polydisc $|w_x|\le R_x$ and hence  $Z_G(\bm w)$ has no zeroes  in  $K_{\bm p}=\prod_{x\in X}[-p_x,p_x]$. Therefore
 $Z_G(\bm w)$ is positive in the whole $K_{\bm p}$ and in particular $Z_G(\bm -\bm p)>0$
and so, by (\ref{p0})  $P(\emptyset)>0$.
On the other hand, for all $S\subset X$, the set $K_{\bm p^S}=\prod_{x\in X}[-p^S_x,p^S_x]\subset K_{\bm p}$. So we also have that
$Z_G(\bm w)$ is positive in $K_{\bm p^S}$ and in particular $Z_G(-\bm p^S)>0$. Therefore,
by (\ref{ps}), $P(S)>0$. Finally, the bound \reff{eq:r2} is a consequence of \reff{eq:r1} and \reff{ps}.
$\Box$

\vv

Estimations on the radius of the analyticity polydisc of $\log Z_G(\bm w)$ is a
classical subject in statistical mechanics.
The best relevant results are contained in the following recent theorem by two of us \cite{FP}.

\begin{theorem}
\label{FP}
Let  ${\bm \mu}=\{\mu_x\}_{x\in {X}}$ be a collection of nonnegative numbers.
Then the function $\log Z_G(\bm w)$ is analytic  for all $\bm w$ such that
\begin{equation}\label{fp}
|w_x|\;\le\;\rfp_x\;\equiv\;\frac{\mu_x}{\ffp_x(\bm \mu)}~~~~~~~~~~~~~~~\forall x\in {X}
\end{equation}
with\
\begin{equation}\label{ffp}
\ffp_x(\bm \mu)=\sum_{R\subset  \Gamma^*_G(x)\atop R\ {\rm indep\ in}\ G}
\prod_{x\in R}\mu_x
\end{equation}
Furthermore,
\begin{equation}
\label{eq:r3}
Z_G(-|\bm w|) \;\ge\; \prod_{x\in X}\left(1-|w_x| \right)^{\ffpt_x(\bm \mu)}
\end{equation}
where
\be\label{fd.0}
\ffpt_x(\bm \mu)=\sum_{R\subset  \Gamma_G(x)\atop R\ {\rm indep\ in}\ G}
\prod_{x\in R}\mu_x
\ee
\end{theorem}

\\{\bf Remark}.
The lower bound (\ref{eq:r3}) is not explicitly given in reference \cite{FP}. It can be proven, however, in a straightforward way from
an upper bound on the (positive) quantity ${-\partial\over\partial |\bm w_i|}\{\ln Z_G(-|\bm w|)\}$ presented in this reference.  For completeness we include the corresponding argument in the appendix.  Note also  that $\ffpt_x(\bm \mu)$ does not depend on $\m_x$ and
\be\label{fft}
\ffp_x(\bm \mu)=\m_x+\ffpt_x(\bm\mu)
\ee

\vv
Plugging Theorem \ref{FP} in Proposition 3, one can replace the Lovasz lemma with the following improved condition.
\begin{theorem}[{\bf Improved Lovasz lemma}]\label{LLL} Suppose that $G$ is a dependence
graph for the family of events $\{A_x\}_{x\in X}$ each one with probability
$\mathbb{P}(A_x)=p_x$ and there exists $\{\mu_x\}_{x \in X}$ real
numbers in $[0,+\infty)$ such that, for each $x\in X$,
\begin{equation}
\label{eq:r4}
              p_x\;\leq\;  \rfp_x\;=\; {\mu_x \over \ffp_x(\bm \mu)}
\end{equation}
where
\begin{equation}\label{ffp11}
\ffp_x(\bm \mu)=\sum_{R\subset  \Gamma^*_G(x)\atop R\ {\rm indep\ in}\ G}
\prod_{x\in R}\mu_x
\end{equation}
Then
\begin{eqnarray}
\label{eq:r3b}
\mathbb{P}(\bar A) 
&\ge & \prod_{x\in X}\left(1-p_x\right)^{\ffpt_x(\bm \mu)}~>~0\label{eq:r3b}
\end{eqnarray}
\end{theorem}

As $\ffp_x(\bm \mu)\le  \varphi_x(\bm \mu)$ [compare (\ref{ab}) and (\ref{ffp11})], $R_x\le R^*_x$.  Thus,
\reff{eq:r4} is less restrictive than the condition  (\ref{lov}) appearing in the Lov\'asz Local Lemma.
Moreover, the later condition depends only on the cardinality
 of the neighborhoods  $\G_G^*(x)$ of each vertex  $x$ of $G$.  In contrast, condition (\ref{eq:r4}) in Theorem \ref{LLL} depends also on
 the edges between vertices of $\G_G^*(x)$.  In consequence, the improvement brought by Theorem \ref{LLL} is maximal
when vertices in each neighborhhod $\G_G^*(x)$ form a clique, and it is null when these vertices form independent sets in the dependency graph.
This occurs, for instance, in bipartite graphs and trees. In particular, if the dependency graph is  a regular tree, 
the Lov\'asz Local Lemma
is tight (see \cite{Sh} or  \cite{Bo} theorem 1.18 or \cite{So} theorem 5.6).

We also point out  that  the lower bound \reff{eq:r3b} is stronger than \reff{a} in the region of validity of the latter.  Indeed, the condition
\begin{equation}\label{eq:rr2}
\frac{1}{1+\mu_x} \;\le\; \left(1-p_x\right)^{\ffpt_x(\bm \mu)}
\end{equation}
holds as as long as
\begin{equation}\label{eq:rr3}
p_x\;\le\;\overline R_x\;\equiv\; 1- \left(\frac{1}{1+\mu_x}\right)^{1/\ffpt_x(\bm \mu)}
\end{equation}
We have
\begin{equation}\label{eq:rr4}
\rfp_x\;\ge\;\overline R_x\;\ge\; \rfpt_x\;\equiv\; {\m_x\over (1+\m_x)\,\ffpt_x(\bm\mu)}\;\ge\; R_x
\end{equation}
The first and second inequalities are a consequence of the elementary identity $(1+a)^b\le 1+ab$ valid for $a\ge -1$ and $0\le b\le 1$.  The last inequality follows easily from the fact that $ (1+\m_x)\,\ffpt_x\ge\ffp_x$.  The inequalities in \reff{eq:rr4} are all strict unless $\ffpt_x(\bm\mu)=1$, i.e. unless $\m_y=0$ for all $y\in \G_G(x)$.
We conclude that the bound \reff{eq:r3b} improves \reff{a} within the region $\{p_x\le R_x\}$ and extends it  to the larger region $\{p_x\le\overline R_x\}$.  The intermediate diameters $\rfpt_x$ could be convenient for calculations.
\medskip

As a final remark we point out that,
 as shown by Scott and Sokal (see Theorem 4.1 in \cite{So}),  Theorem 2 and Proposition 3  also hold
if $G$ is a {\it lopsidependency graph} for the family of events $\{A_x\}_{x\in X}$, i.e.
if, for all $x \in X$ and all $Y \subseteq X\backslash\Gamma^{*}(x)$
we have
$\mathbb{P}(A_x|\bigcap_{y \in Y}\overline{A}_y)\leq \mathbb{P}(A_x)$.  So also the
Lopsided Lov\'{a}sz local lemma can be replaced by the following improved criterion.

\begin{theorem}\label{lop}Let $\{A_x\}_{x\in X}$ be a family of events on some probability
space and $G$ a lopsidependency graph for the family of events with probability
$\mathbb{P}(A_x)=p_x$ for all $x \in X$. Suppose that
$\boldsymbol{\mu}=(\mu_x)_{x\in X}$ are real numbers in
$[0,+\infty)$ such that, for each $i$ we have
\begin{equation}
               p_x\leq {\mu_x \over \ffp_x(\bm \mu)}
\end{equation}
 Then
\begin{equation}\label{eq:r6}
\mathbb{P}(\bar A) \;\ge\; \prod_{x\in X}\left(1-p_x \right)^{\ffpt_x(\bm \mu)} \;>\; 0
\end{equation}
\end{theorem}

\section{Applications}

\\\underline{\it Application 1}.
\vv

\begin{proposition}\label{prop}
Let $G=(V_G,E_G)$ be a graph with maximum degree $\Delta$ and $V_G=V_1
\cup V_2 \cup... \cup V_n$ a partition of vertices set $V_G$ into $n$
pairwise disjoints sets. Suppose that for each set $V_i$ we have
$|V_i| \geq 4\Delta $.
Then, there exists an independent  set $W
\subseteq V_G$ of cardinality $n$ that contains exactly one vertex
from each $V_i$.
\end{proposition}
{\it Proof:} Without loss of generality, we can assume that for all
sets $V_i$ we have $|V_i| = s$ for some integer $s$. The general case follows
from this using the graph induced by $G$ on a union of $n$ subsets
of cardinality $s$, each of them subset of one $V_i$.
We take a random set $W$ of $n$ vertices as follows: we
pick up from each set $V_i$, randomly and independently, a unique
vertex according to a uniform distribution, that is, in
each $V_i$ the probability of a vertex to be chosen is
$1/s$.

Let now  $\tilde E_G=\{e\in E_G: |e\cap V_k|\le 1,  \forall k=1,\dots,n\}$.
Namely, $\Tilde E_G\subset E_G$ contains all edges whose end-points are contained in two distinct  sets $V_i$ and $V_j$ of the partition
of $V_G$.
For  $e=\{a,b\}\in \tilde E_G$, let $W_{e}$ be the event \textquotedblleft $\{a,b\}\subset W$\textquotedblright. Clearly,
$\mathbb{P}(W_{e})=  1/s^{2}\doteq p$ if $e\in \tilde E_G$.  Moreover, if  $e=\{a,b\}\in \tilde E_G$ is such that
$a\in V_i$ and $b\in V_j$, then the event  $W_{e}$ is mutually independent of all the events involving edges whose endpoints
do not lie in does not in $V_i\cup V_j$.

So the graph   $H$ with vertex set $V_H=\tilde E_G$ edge set

$$
E_H=\{\{e,e'\}\subset \tilde E_G:~ |e\cap V_k|+ |e'\cap V_k|> 1~{\rm for\,\, some}~ k=1,\dots,n\}
$$
is a dependency graph
for the family of events
$\{W_{e}\}_{e\in \tilde E_G}$. The maximum degree of the graph $H$ is less than $2 s \Delta$. Indeed, for each of the two vertices of an edge
$e\in \tilde E_G$
there are at most $s$ vertices in the set $V_i$ which contains that vertex, and since
$G$ has degree at most $\D$,  there are at most $\Delta$ edges containing a fixed vertex. It is also clear
that $\Gamma^*_H(e)$ (i.e.  the neighborhood  of $e$ in $H$)   is the union of $2$ subsets, each of cardinality
$s\Delta$ such that  vertices in the same subset are all adjacent.

To apply Theorem \ref{LLL} for the dependence graph $H$,
we take $\mu_e = \mu > 0$, for all $e\in V_H$ and then we
observe that
\begin{eqnarray}\nonumber
\frac{\mu}{\varphi^{*}_e(\mu)} \;\geq\; \frac{\mu}{\left[1+s \Delta \mu\right]^2} \;\equiv\; f(\mu)
\end{eqnarray}
As the right side assumes its maximum value at
$\mu_0=[s\Delta]^{-1}$ we may use theorem \ref{LLL} in the
region
${1\over s^2}\leq {1\over 4s \Delta}=f(\mu_0)
$,
a condition equivalent to say
$
s \geq 4\D
$.
This guarantees that there is a positive  probability that
none of the events $\{W_{e}\}_{e\in \tilde E_G}$ occur. In
other words, we obtain that  there is a positive probability that  $W$
is an independent set that contains one vertex of each $V_i$.
$\qed$\\

We remark  that the same proposition with the $2e$ replacing $4$ is proved in \cite{AS}, page 70, using the Lov\'asz local Lemmma. On the other hand
by a different method it has been shown that the constant can be lowered to $2$  (see \cite{Hax} and references therein).

\vv\vv
\\\underline{\it Application 2: Latin Transversal}.
\vv

Let $A$ be a $n\times n$ matrix with entries $a_{ij}$. Suppose that $a_{ij}$ is integer
for all $i,j=1,\dots, n$. A permutation $\s: \{1,2,\dots, n\}\to \{1,2,\dots, n\}: i\mapsto \s(i)$ is called a Latin transversal of $A$ if the entries
$a_{i\s(i)}$ with $i=1.\dots, n$ are all distinct.
\begin{proposition}
Suppose that $k\le (n-1)/(256/27)$ and suppose that no integer appears in more than $k$ entries of $A$. Then $A$ has a Latin transversal.
\end{proposition}
\\{\it Proof}. Let $\s$ be a permutation of $\{1,2,\dots, n\}$ chosen at random  with uniform distribution. Denote by $T$
the set of all ordered four-tuples $(i,j,i',j')$ such that $i<i'$,  $j\neq j'$ and $a_{ij}=a_{i'j'}$. For each $(i,j,i',j')\in T$, let $A_{iji'j'}$ be the event
that $\s(i)=j$ and $\s(i')=j'$. Clearly $A_{iji'j'}$ has a probability ${1\over n(n-1)}$ to occur and any permutation $\s$ such that $A_{iji'j'}$ occurs is not a Latin transversal. Hence,
a Latin transversal of $A$ exists if there is a non zero probability that none of the events $A_{iji'j'}$ occurs. Following Alon and Spencer \cite{AS}
we define a graph $G$ with vertex set $T$ and two vertices  $(i,j,i',j')$ and  $(p,q,p',q')$ are adjacent if and
only if $\{i, i'\}\cap\{p,p'\}\neq \emptyset$  or $\{j, j'\}\cap\{q,q'\}\neq \emptyset$. This graph has maximum degree less that $4nk$. Indeed, for a fixed
$(i,j,i',j')$, we  can choose $(s,t)$ in $4n$ different ways with
$s\in \{i,i'\}$ or $t\in \{j,j'\}$ for a given $(i,j,i',j')$ and  once $(s,t)$ has been chosen we have less that $k$ choices for $(s',t')$
distinct from $(s,t)$ such that $a_{st}=a_{s't'}$,
since by assumption there are at most $k$ entries of $A$ with the same value. So we have less than $4nk$ four-tuples  $(s,t,s',t')$ such that
$s=i$, or $s=i'$, or $t=j$ or $t=j'$ and such that $a_{st}=a_{s't'}$. Now, to each of those four-tuples  $(s,t,s',t')$ we can associate
uniquely the  four-tuple $(p,q,p',q')= (s,t,s',t')$ if $s<s'$ or the four-tuple $(p,q,p',q')= (s',t',s,t)$ if $s'<s$.

Alon and Spencer showed that $G$ is a lopsidependency graph for the family of events $A_{iji'j'}$. Namely they showed (see \cite{AS} pag. 79 or \cite{ES2}) that
$$
\mathbb{P}\Bigg(A_{iji'j'}|\bigcap_{(p,q,p',q') \in Y}\overline{A}_{pqp'q'}\Bigg)\leq  {1\over n(n-1)}
$$
for any $(i, j , i',j')\in T$  and any set $Y$ of members of $T$ that are nonadjacent in $G$ to
$(i, j , i',j')$.
Therefore we can  apply Theorem \ref{lop} for the lopsidependency graph $G$ with vertex set $T$ described above. We  take  $\mu_x = \mu > 0$, for all $x \in T$ and
we observe that  the set of vertices in $\G^*_G((i,j,i',j'))$ is,  by the previous construction, the union of 4 subsets each of cardinality at most $nk$ and
such that the  all vertices in each one of these four subsets are adjacent.
Thus, for such a graph $G$
$$
\ffp_{(i,j,i',j')}(\mu)\le (1 + nk\mu)^4
$$
and
$$
{\mu\over
\ffp_{(i,j,i',j')}( \mu)}\ge {\mu\over
(1 + nk\mu)^4} \equiv f(\mu)
$$
As the righthand side assumes its maximum value at
$\mu_0=\frac{1}{3nk}$, we may use  Theorem \ref{lop} on the
region $p={1\over n(n-1)}\leq 27/(256nk)=f(\mu_0)$, a condition equivalent to say $k\leq (n-1)/(256/27)$. $\Box$

Note that the same proposition with $4e$ replacing $256/27$ is proved in \cite{AS}, using the Lov\'asz local Lemmma.
\vskip.2cm

\vv\vv
\\\underline{\it Application 3: $k$-SAT}.
Let a set  $U$ of Boolean variables.  Namely, each variable $x\in U$
can take either the value ``true" (indicated by the {\it literal} $x$)  or the value ``false" (indicated by the literal $\bar x$)
Variables or their negations strung together with or symbols are called
{\it clauses}. A set of clauses $C$ joined together by and symbols forms  a Boolean formula  in  {\it Conjunctive Normal Form} (CNF).
A a CNF formula $F=(U, C)$ where $U$ is the set of variables and $C$ is the set of clauses is called a {\it k-SAT} if
all clauses have exactly $k$ variables. A k-SAT is said to be {\it satisfiable} if values true/false can be assigned to its variables in a way that makes the formula logically true
(i.e. every clause in $C$ contains at least one literals taking the value``true").

\begin{theorem}
If F = (U,C) is an instance of k-SAT such that each variable lies in at most ${1\over e}
2^k/k$ clauses,
then F is satisfiable.
\vv
\end{theorem}

\\{\bf Proof}. Suppose that each variables $x\in U$  in a clause $c\in C$ is true with probability $1/2$. For each clause $c\in C$
let $A_c$ be the event that $c$ is false (i.e. all the $k$ variables in $c$ are set at the value ``false"). Suppose now that
each variables $x\in U$ is contained in at most $x$ clauses. Then the event $A_c$ is independent of all other events involving clauses
not containing the variables of $c$. Hence a dependency graph for the set of events $\{A_c\}_{c\in C}$ is the graph $H=(C,E)$ with vertex
set $C$ (the set of clauses of $F$)
and edge set
$$
E=\{\{c,c'\}\subset C: \mbox{$c$ and $c'$ have at least one variable in common}\}
$$
We have clearly that $Prob(A_c)=2^{-k}$ and the theorem is proven if we can prove that there is a positive probability that
none of the events $\{A_c\}_{c\in C}$  occur. We apply Theorem \ref{LLL}.
Since by hypothesis each variables $x\in U$ is contained in at most $N$ clauses and each clause contains $k$ variables, it is clear
that the neighbor $\G^*_H(c)$ of a vertex $c$ of the dependency graph $H$ is formed  by at most $k$ cliques each one containing at most $N$ vertices.
Hence choosing $\m_c=\m$ for all $c\in C$,  condition \reff{eq:r4} Theorem \ref{LLL} is satisfied if
\begin{equation}\label{ksat}
2^{-k} \le {\mu\over (1+ N\m)^k }
\end{equation}
The maximum of r.h.s. of inequality \reff{ksat} is attained for $\m={1/(k-1)n}$ and substituting this value inequality \reff{ksat}
becomes
$$
2^{-k} \le {1\over (k-1)N(1+ {1\over k-1})^k }=  {1\over kN(1+ {1\over k-1})^{k-1} }
$$
which is satisfied if  $N\le {2^k\over e k}$. $\Box$

Note that the same proposition with $4$ replacing $e$ is proved in \cite{Mo}   (see there theorem 3.1) , using the Lov\'asz local Lemmma.

\renewcommand{\theequation}{A.\arabic{equation}}
\setcounter{equation}{0}  
\section*{Appendix. Proof of bound (\ref{eq:r3b})}  

The partition functions of hard-core gases  satisfy, in the sense of formal power series, the identities (see e.g. \cite{Ru}, \cite{FP})
$$
-{\partial\over \partial \r_{x}}\log Z_G(-\bm \r)= \Pi_x(\bm \r)\qquad\qquad \forall\,x\in X
$$
$\Pi_x(\bm \r)$ are ($G$-dependent) formal power series on $\bm \r=\{\r_x\}_{x\in X}$, all whose terms are positive if $\r_y\ge 0$, $y\in X$. The simultaneous convergence for  some $\bm \r$ of all the series  $\Pi_y(\bm \r)$, $y\in X$, implies the analyticity of $\log Z_G(\bm w)$   in the polydisc  $\{|w_x|\le \r_x \}_{x\in X}$.  Furthermore,  within the region of convergence,
\begin{equation}\label{mlog0}
|\log Z_G(\bm w)|\;\le\; -\log Z_G(-|\bm w|)
\end{equation}
 and
\begin{equation}
-\log Z_G(-\bm \r)+ \log Z_{G\backslash\{x\}}(-\bm \r) \;=\; \Sigma_x(\bm \r)
\end{equation}
where $G\backslash\{x\}$ denotes the restriction of $G$ to $X\backslash \{x\}$ and
\begin{equation}\label{sigma}
\Sigma_x(\bm \r)\;=\;
\r_{x}\int_0^1 \Pi_{x}(\bm \r(\a))d\a
\end{equation}
with
$$
\r_y(\a)\;=\;
\begin{cases}\r_y &\text{if $y\neq x$}\cr
\a \r_y &\text{if $y=x$}
\end{cases}
$$
As a consequence
\be\label{mloga}
-\log Z_G(-|\bm w|)\;\le\; \sum_{x\in X} \Sigma_x(\bm \r)
\ee
The lower bound \reff{eq:r3b} is proven from upper bounds for $\Sigma_x(\bm \r)$.
They follow from the following proposition proved in \cite{FP}.

\begin{proposition}
\label{fixedp}
Let $\bm\mu\equiv \{\mu_x\ge 0\}_{x\in X}$ and let, for all $x\in X$, $\r_x\le \rfp$ (with $\rfp$ defined in (\ref{FP})).
Then there exists a convergent positive term series  $\pfp_x(\bm \r)$  monotonic in $\bm \r$
such that
\begin{enumerate}
\item[i)]
If $\bm \r\bm\Pi^{^{_{\rm FP}}}\!(\bm\r)=\{\r_x\pfp_x\!(\bm \r)\}_{x\in X}$ and $\fbfp\!(\bm \mu)= \{\ffp_x\!(\bm \mu)\}_{x\in X}$, then
$\bm \r \pbfp\!(\bm\r)$
is fixed point of the map
\begin{eqnarray}
\label{eq:rr7}
\bm T^{\bm \r}\equiv \{T_x^{\bm \r}\}_x:  [0,\infty)^X&\longrightarrow& [0,\infty)^X\nonumber\\
T_x^{\bm \r}\!(\bm\cdot)&=&\r_x \ffp_x\!(\bm \cdot)
\end{eqnarray}
 i.e.
\begin{equation}\label{fixed}
\r_x\pfp_x(\bm \r)= T^{\bm\r}_{x}\Big(\bm \r\bm\pbfp(\bm\r)\Big)
\end{equation}
\item[ii)] The following bounds hold
$$
\r_x\Pi_x(\bm \r)\;\le\; \r_x\pfp_x\!(\bm \r)\;\le\; \mu_x
$$
\end{enumerate}
\end{proposition}

\\\\Item {\it ii)} of Proposition \ref{fixedp} immediately implies that $\log Z_G(\bm w)$ is analytic in the polydisk  $\{|w_x|\le R_x\}_{x\in X}$.
\\To obtain an upper bound for $\Sigma_x(\bm \r)$  we make use of item
{\it i)}. Indeed, by (\ref{fixed}) and (\ref{fft})
\begin{eqnarray}
\pfp_x\!(\bm \r) ~=~ \ffp_x\!(\bm \r \pfp\!(\bm \r))
&= & \r_x \pfp_x\!(\bm \r) +\ffpt_x\!(\bm \r \pbfp\!(\bm \r))\nonumber\\
&\le & \r_x \pfp_x\!(\bm \r)+
\ffpt_x\!(\bm \mu)\nonumber
\end{eqnarray}
Thus, since $\Pi_x(\bm \r)\le \pfp_x\!(\bm \r)$,
$$
{\P_x(\bm \r)}\;\le\;  {\ffpt_x\!(\bm \mu)\over 1- \r_x}
$$
By equation (\ref{sigma}) this implies
$$
\Sigma_x (\bm \r)\;\le\;\r_x\int_0^1 {\ffpt_x\!(\bm \mu)\over 1- \a\r_x}d\a\;=\;-\ln \left[1- \r_x\right]^{\ffpt_x\!(\bm \mu)}
$$
which, by (\ref{mloga}), implies that for all $\r_x\le \rfp_x$
$$
Z_G(-\bm \r) \ge \prod_{x\in X}\left[1- \r_x\right]^{\ffpt_x\!(\bm \mu)}
$$
\section*{Acknowledgments}
It is a pleasure to thank Yoshiharu Kohayakawa for helpful comments.
This work  has been partially supported by
Conselho Nacional de Desenvolvimento Cient\'{i}fico e Tecnol\'ogico (CNPq),
 CAPES (Coordena\c{c}\~ao de Aperfei\c{c}oamento de Pessoal de
N\'{\i}vel Superior, Brasil) and FAPEMIG (Funda{c}\~ao de Amparo \`a Pesquisa do Estado de Minas Gerais). RB
thanks the hospitality of Laboratoire de Math\'{e}matiques
Rapha\"{e}l Salem in the Universit\'{e} de Rouen during his PhD work.
RF and AP  thank The Newton Institute for the generous support
during the Combinatorics and Statistical Mechanics programme  and the Mathematics Department
of Universit\`a di Roma Tor Vergata (AP).

\end{document}